\documentclass[11pt]{article}
\usepackage{amssymb}
\usepackage{amsmath}
\usepackage{amsfonts}

\begin{document}

\title{ \textbf{Branching rules for }$S_{2N}\rightarrow W\left( B_{N}\right)
$}
\author{Godofredo Iommi Amun\'{a}tegui \\
{\small \emph{Instituto de F\'{\i}sica, P. Universidad Cat\'{o}lica de
Valpara\'{\i}so, Casilla 4059, Valpara\'{\i}so, Chile.}} \\
\texttt{godofredo.iommi@pucv.cl}}
\date{}
\maketitle

\begin{abstract}
This note presents a procedure to determine the reduction of the irreducible
and the induced characters of the symmetric group $S_{2N}$ in terms of the
irreducible and induced characters of the hyperoctahedral group $W\left(
B_{N}\right) =Z_{2}^{N}\sim S_{N}$.

\noindent Mathematical Subject Classification

\noindent $20B_{xx},$ $20C_{xx},$ $20E_{xx}$

\end{abstract}

\section{ Introduction}

To each classical Lie group corresponds a finite group generated by the
reflections of its root system, called the Weyl group. There has been a
number of situations in which the Weyl groups have played an important role.
This importance grew out of the various possibilities of application to
physical problems i.e., particle physics, discrete $\sigma $ models, lattice
gauge theories, chiral models (Ref. \cite{Mandula,Baake}).

The Symmetric group $S_{N}$ is the Weyl group of the Unitary Group. For $%
B_{N}=SO\left( 2N+1\right) $ and $C_{N}=S_{p}\left( 2N\right) $, the Weyl
Groups $W\left( B_{N}\right) $ and $W\left( C_{N}\right) $ are isomorphic. $%
W\left( B_{N}\right) $\ is $Z_{2}^{N}\sim S_{N}$, the wreath product of the
abelian group $Z_{2}^{N}$ generated by the $N$ sign changes $\left(
+i,-i\right) $, $1\leq i\leq N$, and the symmetric group $S_{N}$. The order
of \ $W\left( B_{N}\right) $\ is $2^{N}N!$\ (Ref. \cite{Humphreys,James}).
Let $K_{N}$\ be defined as the convex hull of points $\pm e_{i}$, $1\leq
i\leq N$, where $e_{1},\ldots e_{N}$\ are the unit coordinate vectors in $%
R^{N}$. It is the N-dimensional generalization of the octahedron $K_{3}$.
The group of symmetries of $K_{N}$, called the hyperoctahedral group is $%
W\left( B_{N}\right) $. The structure and representation of this group have
been studied (Ref. \cite{Stanley,Keber,Geissinger}). Moreover the
hyperoctahedral groups appear in numerous applications such as weakly bound
water clusters, non-rigid molecules, disordered proteins and the enumeration
of isomers (Ref. \cite{Balasubramanian1,Balasubramanian2}). The
hyperoctahedral group $Z_{2}^{N}\sim S_{N}$\ is a subgroup of the symmetric
group $S_{2N}$. The purpose of this note is to propose a procedure to solve
the reduction $S_{2N}\rightarrow \left( Z_{2}^{N}\sim S_{N}\right) $.
Although there are already computer codes available to generate the
character tables of $S_{N}$\ for any $N$, and their wreath products (Ref.
\cite{GAP}), to the best of my knowledge this branching case has not been
treated as yet.

In order to make this article reasonably self-contained some pertinent
results already published will be exposed anew. In Section 2 and Section 3,
respectively, algorithms for the irreducible and induced characters of $%
S_{2N}$ and $W\left( B_{N}\right) $ are treated. Section 4 deals with the
reduction $S_{2N}\rightarrow W\left( B_{N}\right) $.

\section{The induced and the irreducible characters}

Consider a partition $\left( \lambda \right) =\left( \lambda _{1},\ldots
,\lambda _{p}\right) $\ of $\ 2N$, where $\lambda _{1}+\lambda _{2}+\ldots
+\lambda _{p}=2N$, $\lambda _{1}\geq \lambda _{2}\geq \ldots \geq \lambda
_{p}=0$; $p\left( 2N\right) $\ is the number of partitions of $2N$.

Corresponding to each partition of $2N$\ we can construct $S_{\lambda
_{1}}\times S_{\lambda _{2}}\times \ldots \times S_{\lambda _{p}}$. Such
subgroups are called the canonical subgroups of $S_{N}$. Let $C$ be a class
of $S_{2N}$\ characterized by its cycle structure $\left( 1^{\alpha
},2^{\beta },3^{\gamma },\ldots \right) $. This symbol denotes that the
permutations in $C$ contain $\alpha $\ 1-cycles, $\beta $\ 2-cycles, $\gamma
$ 3-cycles, etc., where $\alpha +2\beta +3\gamma +\ldots =2N$. Besides for
each $S_{\lambda _{1}}$ we have%
\begin{equation}
\alpha _{i}+2\beta _{i}+3\gamma _{i}+\ldots =\lambda _{i}  \tag{A1}
\label{A1}
\end{equation}

The character induced in $S_{2N}$ by the identity representation of a
canonical subgroup is%
\begin{equation*}
\phi _{\left( 1^{\alpha },2^{\beta },\ldots \right) }^{\left( \lambda
\right) }=\sum \frac{\alpha !}{\alpha _{1}!\alpha _{2}!\ldots }\cdots \frac{%
\beta !}{\beta _{1}!\beta _{2}!\ldots }\cdots \frac{\gamma !}{\gamma
_{1}!\gamma _{2}!\ldots }\cdots
\end{equation*}

Where
\begin{equation}
\sum \alpha _{i}=\alpha \,,\qquad \sum \beta _{i}=\beta \,,\qquad \sum
\gamma _{i}=\gamma \,,\ldots  \tag{A2}  \label{A2}
\end{equation}

The sum is over all the integer solutions of the system of Eqs. (\ref{A1})
and (\ref{A2}). These characters may be arranged as the entries of a $%
p\left( 2N\right) \times p\left( 2N\right) $\ matrix $\phi $ whose rows and
columns are labeled, respectively, by partitions of $2N$ arranged in
lexicographical order and by the classes (Ref. \cite{Iommi}).

The table of irreducible characters of $S_{2N}$\ may be derived from $\phi $
(Ref. \cite{Iommi}). Each row $\phi _{i}$ must be considered as a vector; it
suffices to orthonormalize them via the Gram-Schmidt method to get the rows $%
x_{i}$ of the irreducible characters table X, i.e.,%
\begin{equation}
x_{i}=\phi _{i}-\overset{i-1}{\underset{k=1}{\sum }}\left( \phi
_{i}Kx_{k}\right) x_{k}  \label{uno}
\end{equation}

(for $i=1,$ $x_{i}=\phi _{1}$), where $x_{i}$ and $\phi _{i}$ are the i-th
rows of $X$ and $\phi $ respectively, and $K$ is a diagonal matrix whose
elements are%
\begin{equation*}
\left[ K_{ij}\right] =\delta _{jk}\,\frac{C}{\left( 2N\right) !}
\end{equation*}%
$C$ is the order of the class $\left( 1^{\alpha },2^{\beta },3^{\gamma
},\ldots \right) $ of $S_{2N}$, $C=\frac{\left( 2N\right) !}{1^{\alpha
}\alpha !,2^{\beta }\alpha !\ldots ,}.$

Expression (\ref{uno}) may be written as%
\begin{equation}
\phi _{i}=x_{i}+\overset{i-1}{\underset{k=1}{\sum }}\left( \phi
_{i}Kx_{k}\right) x_{k}
\end{equation}

Considering the coefficients of the $x_{k}$ we get a lower triangular matrix
$\Delta $ such that det $\Delta =1$. In general we have for $S_{2N}$%
\begin{equation}
\phi =\Delta X
\end{equation}

As an example for $S_{4}$ we have:

\bigskip

\begin{tabular}{|l|l|l|l|l|lll|l|l|l|l|l|lll|l|l|l|l|l|}
\cline{1-5}\cline{5-5}\cline{9-9}\cline{9-13}\cline{12-13}\cline{17-21}
1 & 1 & 1 & 1 & 1 &  &  &  & 1 & 0 & 0 & 0 & 0 &  &  &  & 1 & 1 & 1 & 1 & 1
\\ \cline{1-5}\cline{5-5}\cline{9-9}\cline{9-13}\cline{12-13}\cline{17-21}
4 & 2 & 0 & 1 & 0 &  &  &  & 1 & 1 & 0 & 0 & 0 &  &  &  & 3 & 1 & -1 & 0 & -1
\\ \cline{1-5}\cline{5-5}\cline{9-9}\cline{9-13}\cline{12-13}\cline{17-21}
6 & 2 & 2 & 0 & 0 &  & = &  & 1 & 1 & 1 & 0 & 0 &  &  &  & 2 & 0 & 2 & -1 & 0
\\ \cline{1-5}\cline{5-5}\cline{9-9}\cline{9-13}\cline{12-13}\cline{17-21}
12 & 2 & 0 & 0 & 0 &  &  &  & 1 & 2 & 1 & 1 & 0 &  &  &  & 3 & -1 & -1 & 0 &
1 \\ \cline{1-5}\cline{5-5}\cline{9-9}\cline{9-13}\cline{12-13}\cline{17-21}
24 & 0 & 0 & 0 & 0 &  &  &  & 1 & 3 & 2 & 3 & 1 &  &  &  & 1 & -1 & 1 & 1 & 1
\\ \cline{1-5}\cline{5-5}\cline{9-9}\cline{9-13}\cline{12-13}\cline{17-21}
\end{tabular}

\bigskip

\section{The Induced and the Irreducible characters of $W\left( B_{N}\right)
$}

The set of all $g=\left( \sigma ;f\right) $, where $\sigma \in S_{2N}$ and $%
f $ is a mapping of $\left[ 1,2N\right] $ into $Z_{2},$ together with the
composition defined by%
\begin{equation*}
\left( \sigma ^{\prime };f^{\prime }\right) \left( \sigma ;f\right) =\left(
\sigma ^{\prime }\sigma ;f^{\prime }\left( f\sigma ^{\prime ^{-1}}\right)
\right)
\end{equation*}%
form the group $W\left( B_{N}\right) =Z_{z}^{N}\sim S_{N}.$

The cycles of the permutation are called \textquotedblleft cycles of $g$%
\textquotedblright . A cycle $\left( a_{1,}\ldots ,a_{\beta }\right) $\ of $%
g $ is positive or negative if $f\left( a_{1}\right) \ldots f\left( a_{\beta
}\right) =+1$ or $-1$. Let $\beta =\left( \beta _{1},\ldots ,\beta
_{k}\right) $\ be the $\beta $ system of cycles of $\sigma $, and suppose
the cycles are arranged in such a way that a negative cycle necessarily
precedes a positive cycle of equal length. Then $\left( \beta ,b\right) $ is
called the $\beta $ system of cycles of $g$, where $b:=\left( b_{1},\ldots
,b_{k}\right) $ with $b_{i}:=1$\ or 0 if the i-th cycle is positive or
negative (remark: if $\beta _{i}=\beta _{:+1}$, then $b_{i}\leq b_{i+1}$).
Moreover if $\alpha _{i}^{+}$ and $\alpha _{i}^{-}$ denote then number of
positive and negative cycles, respectively, of length $i$ of $g$, then

\begin{equation*}
\alpha =\left( \alpha _{1}^{+},\alpha _{1}^{-},\alpha _{2}^{+},\alpha
_{2}^{-},\ldots ,\alpha _{\ell }^{+},\alpha _{\ell }^{-}\right)
\end{equation*}

is called the $\alpha $ system of cycles of $g$ (remark: if $\alpha
_{i}:=\alpha _{i}^{+}+\alpha _{i}^{-}$ then $\underset{i}{\sum }i\alpha
_{i}^{+}=N$).

The elements of $W\left( B_{N}\right) $ are conjugates $iff$ the have the
same $\alpha $ system of cycles and $iff$ they have the same $\beta $ system
of cycles. The class of elements with $\alpha $ system $\alpha =\left(
\alpha _{1}^{+},\ldots ,\alpha _{1}^{-}\right) $ is denoted $C(\alpha )$.

Let $\lambda =\left( \lambda _{1},\ldots ,\lambda _{k}\right) \ $be a
partition of $N\left( \lambda _{1}\geq \lambda _{2}\geq \ldots \geq \lambda
_{k}\right) $ and $b=\left( b_{1},\ldots ,b_{k}\right) $ be such that $%
b_{i}=1$ or $0$ (remark: if $\lambda _{i}=\lambda _{i+1}$, then $b_{i}\leq
b_{i+1}$). The subgroup $\left( Z_{2}^{\left( \lambda _{1}-b_{1}\right)
}\sim S_{\lambda _{1}}\right) \times \left( Z_{2}^{\left( \lambda
_{2}-b_{2}\right) }\sim S_{\lambda _{1}}\right) \ldots $, denoted by $%
S\left( \lambda ,b\right) $, is a canonical subgroup of $W\left(
B_{N}\right) $. Then, for the class $C(\alpha )$ and the canonical subgroup $%
S\left( \lambda ,b\right) $ the algorithm giving the character $I_{S\left(
\lambda ,b\right) }^{\left( C\left( \alpha \right) \right) }$ of the
representation of $W\left( B_{N}\right) $ induced by the identity
representation of $S\left( \lambda ,b\right) $ is:

\begin{equation*}
I_{S\left( \lambda ,b\right) }^{\left( C\left( \alpha \right) \right)
}=2^{\left( \underset{i}{\Sigma }b_{i}\right) }\left( \sum \frac{\Pi
_{i=1}^{\ell }\left( \alpha _{i}^{+}\right) !\left( \alpha _{i}^{-}\right) !%
}{\Pi _{i=1}^{\ell }\Pi _{j=1}^{k}\left( \alpha _{ij}^{+}\right) !\left(
\alpha _{ij}^{-}\right) !}\right)
\end{equation*}

Then sum concerns the matrices $\left( \alpha _{ij}^{+/-}\right) $ of dim $%
\ell \times k\times 2$ where%
\begin{equation*}
\forall i_{0},\text{ \ }\underset{j=1}{\overset{k}{\sum }}\alpha
_{i_{0}j}^{+}=\alpha _{i_{0}}^{+}
\end{equation*}%
and%
\begin{eqnarray*}
\underset{j=1}{\overset{k}{\sum }}\alpha _{i_{0}j}^{-} &=&\alpha _{i_{0}}^{-}
\\
\forall j_{0}\,,\text{ \ }\underset{j=1}{\overset{\ell }{\sum }}i\left(
\alpha _{ij_{0}}^{+}+\alpha _{ij_{0}}^{-}\right) &=&\lambda j_{0}
\end{eqnarray*}%
Besides $\forall j_{0},$ if $bj_{0}=1,$ then $\underset{i}{\sum }\alpha
_{ij_{0}}^{-}$ is an even number. The order of the class $C\left( \alpha
\right) $ is%
\begin{equation*}
\left\vert C\left( \alpha \right) \right\vert =N!\overset{\ell }{\underset{%
i=1}{\prod }}\left( \frac{2^{\alpha _{i}\left( i-1\right) }}{i^{\alpha
_{i}}\left( \alpha _{i}^{+}\right) !\left( \alpha _{i}^{-}\right) !}\right)
\end{equation*}

By means of such an algorithm, the induced character table $I\left\{ W\left(
B_{N}\right) \right\} $ is obtained. Each row of the table is given by the
corresponding $I_{S\left( \lambda ,b\right) }\left( C\left( \alpha \right)
\right) .$

For $N=2$, the table of induced characters is:

\bigskip

\begin{tabular}{llclllll}
&  & 1 & 2 & 1 & 2 & 2 & classes order \\ \cline{3-3}\cline{3-7}
&  & \multicolumn{1}{|c}{1} & \multicolumn{1}{|l}{1} & \multicolumn{1}{|l}{1}
& \multicolumn{1}{|l}{1} & \multicolumn{1}{|l|}{1} & \multicolumn{1}{|l}{}
\\ \cline{3-3}\cline{3-7}
&  & \multicolumn{1}{|c}{2} & \multicolumn{1}{|l}{0} & \multicolumn{1}{|l}{2}
& \multicolumn{1}{|l}{2} & \multicolumn{1}{|l|}{0} & \multicolumn{1}{|l}{}
\\ \cline{3-3}\cline{3-7}
$I\left\{ W\left( B_{2}\right) \right\} $ & = & \multicolumn{1}{|c}{2} &
\multicolumn{1}{|l}{2} & \multicolumn{1}{|l}{2} & \multicolumn{1}{|l}{0} &
\multicolumn{1}{|l|}{0} & \multicolumn{1}{|l}{} \\ \cline{3-3}\cline{3-7}
&  & \multicolumn{1}{|c}{4} & \multicolumn{1}{|l}{2} & \multicolumn{1}{|l}{0}
& \multicolumn{1}{|l}{0} & \multicolumn{1}{|l|}{0} & \multicolumn{1}{|l}{}
\\ \cline{3-3}\cline{3-7}
&  & \multicolumn{1}{|c}{8} & \multicolumn{1}{|l}{0} & \multicolumn{1}{|l}{0}
& \multicolumn{1}{|l}{0} & \multicolumn{1}{|l|}{0} & \multicolumn{1}{|l}{}
\\ \cline{3-3}\cline{3-7}
\end{tabular}

\bigskip

The table of irreducible characters $Y\left\{ W\left( B_{N}\right) \right\} $
can be obtained from $I\left\{ W\left( B_{N}\right) \right\} $. As before
each row of $I\left\{ W\left( B_{N}\right) \right\} $ most be considered as
a vector and via the Gram-Schmidt procedure the rows of $Y\left\{ W\left(
B_{N}\right) \right\} $ are obtained. In general

\begin{equation*}
Y_{i}=I_{i}-\underset{k=1}{\overset{i-1}{\sum }}\left( I_{i}DY_{k}\right)
Y_{k}\text{ \ \ \ \ \ \ \ \ \ \ \ \ \ \ \ \ \ \ }\left( \text{for \ }i=1,%
\text{ \ \ }Y_{1}=I_{1}\right)
\end{equation*}%
where $Y_{i}$ and $I_{i}$ are the i-th row of $Y\left\{ W\left( B_{N}\right)
\right\} $ respectively and $D$ is a triangular matrix whose elements are $%
\left( D_{\alpha \beta }\right) =\delta _{\alpha \beta }\frac{\left\vert
C\left( \alpha \right) \right\vert }{2^{N}N!},$ $\left\vert C\left( \alpha
\right) \right\vert $ is the order of the class $C\left( \alpha \right) $ of
$W\left( B_{N}\right) .$ Here

\begin{equation}
I\left( W\left( B_{N}\right) \right) =DY\left\{ W\left( B_{N}\right)
\right\} .
\end{equation}

For instance, for $W\left( B_{2}\right) $, the Weyl group of $SO\left(
5\right) $, we have\bigskip

\begin{equation*}
\begin{tabular}{|l|l|l|l|l|lll|l|l|l|l|l|l|l|l|l|l|l|}
\cline{1-5}\cline{9-13}\cline{15-19}
1 & 1 & 1 & 1 & 1 &  &  &  & 1 &  &  &  &  &  & 1 & 1 & 1 & 1 & 1 \\
\cline{1-5}\cline{9-13}\cline{15-19}
2 & 0 & 2 & 2 &  &  &  &  & 1 & 1 &  &  &  &  & 1 & -1 & 1 & 1 & -1 \\
\cline{1-5}\cline{9-13}\cline{15-19}
2 & 2 & 2 &  &  &  & = &  & 1 & 0 & 1 &  &  &  & 1 & 1 & 1 & -1 & -1 \\
\cline{1-5}\cline{9-13}\cline{15-19}
4 & 2 &  &  &  &  &  &  & 1 & 0 & 1 & 1 &  &  & 2 & 0 & -2 & 0 & 0 \\
\cline{1-5}\cline{9-13}\cline{15-19}
8 &  &  &  &  &  &  &  & 1 & 1 & 1 & 2 & 1 &  & 1 & -1 & 1 & -1 & -1 \\
\cline{1-5}\cline{9-13}\cline{15-19}
\end{tabular}%
\end{equation*}

\bigskip

i.e., $I\left\{ W\left( B_{2}\right) \right\} =DY\left\{ W\left(
B_{2}\right) \right\} .$

\section{The Reduction $S_{2N}\rightarrow W\left( B_{N}\right) $}

In this section we expose a procedure to express the content of the
irreducible and the induced characters of $S_{2N}$ in terms of the
irreducible and the induced characters of its subgroup $W\left( B_{N}\right)
$. Such an algorithmic process is valid in general i.e., for any $N$.
However it must be pointed out that every branching case must be treated
with due regard to its own structural traits (see Appendix (I)).

As a matter of fact we shall envisage the reduction from two different
points of wiew (hereafter Method (A) and Method (B)).

\subsection{Method (A)}

We already know that for $S_{2N}$ we have $\phi =\Delta X$ (Section 2) and
for $W\left( B_{N}\right) $ $I=DY$ (Section 3). Besides, in order to carry
out the reduction use must be made of the modified characters tables $%
X^{\prime }$ and $\phi ^{\prime }$ (se Appendix (I)). The characters of $%
S_{2N}$ can be expressed in terms of the characters of $W\left( B_{N}\right)
$ by means of reduction matrices. We denote the reduction matrices for the
irreducible and induced characters as

\begin{equation*}
R_{Y_{W\left( B_{N}\right) }}^{X_{S_{2N}}}\text{ \ \ \ and \ \ \ }%
R_{I_{W\left( B_{N}\right) }}^{\phi _{S_{2N}}}
\end{equation*}

\ \ \ \ \ (in short, $R_{1}$ and $R_{2}$ respectively) Then:

\begin{equation}
X^{\prime }=R_{1}Y
\end{equation}

\begin{equation}
\phi ^{\prime }=R_{2}I  \label{seis}
\end{equation}

To obtain the entries of the reduction matrices a system of $P\left(
N\right) $ linear equations with $K\left( W\left( B_{N}\right) \right) $
unknowns must be solved via $K\left( W\left( B_{N}\right) \right) $
independent linear equations.

$K\left( W\left( B_{N}\right) \right) $ is the number of classes of $W\left(
B_{N}\right) $. A simple expression for $K\left( W\left( B_{N}\right)
\right) $ appears in ref \cite{Balasubramanian2}.

Let us note that (\ref{seis}) can be written as

\begin{equation}
\phi ^{\prime }=R_{2}I=R_{2}DY
\end{equation}

and

\begin{equation*}
\phi ^{\prime }=\Delta ^{\prime }X^{\prime }=\Delta ^{\prime }R_{1}Y
\end{equation*}

then

\begin{equation*}
R_{2}DY=\Delta ^{\prime }R_{1}Y
\end{equation*}

hence

\begin{equation}
R_{2}D=\Delta ^{\prime }R_{1}  \label{ocho}
\end{equation}

This equation establishes a direct relation between the two branching
matrices.

To illustrate equation (\ref{ocho}), we shall consider the simplest
reduction case $S_{4}\rightarrow W\left( B_{2}\right) :$

\bigskip

\begin{tabular}{|l|l|l|l|l|l|l|l|l|l|l|l|l|l|l|l|l|l|l|l|l|l|l|}
\cline{1-5}\cline{5-5}\cline{7-7}\cline{7-11}\cline{10-11}\cline{13-17}\cline{13-17}\cline{19-23}
1 &  &  &  &  &  & 1 &  &  &  &  &  & 1 &  &  &  &  &  & 1 &  &  &  &  \\
\cline{1-5}\cline{5-5}\cline{7-7}\cline{7-11}\cline{10-11}\cline{13-17}\cline{13-17}\cline{19-23}
&  &  & 1 &  &  & 1 & 1 &  &  &  &  & 1 & 1 &  &  &  &  &  &  & 1 & 1 &  \\
\cline{1-5}\cline{5-5}\cline{7-7}\cline{7-11}\cline{10-11}\cline{13-17}\cline{13-17}\cline{19-23}
& 1 &  & 1 &  &  & 1 & 0 & 1 &  &  & = & 1 & 1 & 1 &  &  &  & 1 & 1 &  &  &
\\
\cline{1-5}\cline{5-5}\cline{7-7}\cline{7-11}\cline{10-11}\cline{13-17}\cline{13-17}\cline{19-23}
&  &  & 1 & 1 &  & 1 & 0 & 1 & 1 &  &  & 1 & 2 & 1 & 1 &  &  &  &  &  & 1 & 1
\\
\cline{1-5}\cline{5-5}\cline{7-7}\cline{7-11}\cline{10-11}\cline{13-17}\cline{13-17}\cline{19-23}
&  &  &  & 3 &  & 1 & 1 & 1 & 2 & 1 &  & 1 & 3 & 2 & 3 & 1 &  &  & 1 &  &  &
\\
\cline{1-5}\cline{5-5}\cline{7-7}\cline{7-11}\cline{10-11}\cline{13-17}\cline{13-17}\cline{19-23}
\end{tabular}

\bigskip

\subsection{Method (B)}

This approach relies on two branching rules which have been solved. The
first one is the classic Weyl's rule for $S_{N}\rightarrow S_{N-1}$: "The
irreducible representations of $S_{N}$ with the symmetry pattern $\left(
\lambda _{1},\lambda _{2},\lambda _{3},\ldots \right) $ reduces on
restricting $S_{N}$ to the subgroup $S_{N-1}$ associated with the patterns $%
\left( \lambda _{1}-1,\lambda _{2},\lambda _{3},\ldots \right) ;$ $\left(
\lambda _{1},\lambda _{2}-1,\lambda _{3},\ldots \right) ;$ $\left( \lambda
_{1},\lambda _{2},\lambda _{3}-1,\ldots \right) $ and so on. Those patterns
in which the rows are not arranged in decreasing lenght are to be omitted"
(Ref. \cite{Weyl}). Such a reduction may be written as a matrix whose rows
and columns are indexed by the partitions of $N$ and $N-1$ ordered in
lexicographic orden. For example the matrix corresponding to $%
S_{4}\rightarrow S_{3}$ is:

\bigskip

\begin{tabular}{llll}
& 3 & 21 & 111 \\ \cline{2-4}
\multicolumn{1}{r}{4} & \multicolumn{1}{|l}{1} & \multicolumn{1}{|l}{} &
\multicolumn{1}{|l|}{} \\ \cline{2-4}
\multicolumn{1}{r}{31} & \multicolumn{1}{|l}{1} & \multicolumn{1}{|l}{1} &
\multicolumn{1}{|l|}{} \\ \cline{2-4}
\multicolumn{1}{r}{22} & \multicolumn{1}{|l}{} & \multicolumn{1}{|l}{1} &
\multicolumn{1}{|l|}{} \\ \cline{2-4}
\multicolumn{1}{r}{211} & \multicolumn{1}{|l}{} & \multicolumn{1}{|l}{1} &
\multicolumn{1}{|l|}{1} \\ \cline{2-4}
\multicolumn{1}{r}{1111} & \multicolumn{1}{|l}{} & \multicolumn{1}{|l}{} &
\multicolumn{1}{|l|}{1} \\ \cline{2-4}
\end{tabular}

\bigskip

\bigskip

The second one is the reduction rule for the hyperoctahedral group (Ref.
\cite{Doeraene}). We have then :\bigskip

(a)\qquad $S_{N}\rightarrow S_{N-1}\rightarrow \ldots \rightarrow
S_{2}\rightarrow S_{1}$

\bigskip

(b)\qquad $W\left( B_{N}\right) \rightarrow W\left( B_{N-1}\right)
\rightarrow \ldots \rightarrow W\left( B_{1}\right) $\bigskip

Since $S_{2}$ and $W\left( B_{1}\right) $ are isomorphic, from (a) and (b)
we deduce

\begin{equation*}
\left\{ S_{2N}\rightarrow W\left( B_{N}\right) \right\} \left\{ W\left(
B_{N}\right) \rightarrow W\left( B_{1}\right) \right\} =S_{2N}\rightarrow
W\left( B_{N}\right)
\end{equation*}

For $N=2$%
\begin{equation*}
\left\{ S_{4}\rightarrow W\left( B_{2}\right) \right\} \left\{ W\left(
B_{2}\right) \rightarrow W\left( B_{1}\right) \right\} =S_{4}\rightarrow
W\left( B_{1}\right)
\end{equation*}

\bigskip

(i)\qquad $W\left( B_{2}\right) \rightarrow W\left( B_{1}\right) $

\qquad \qquad
\begin{tabular}{|l|l|}
\hline
1 &  \\ \hline
& 1 \\ \hline
1 &  \\ \hline
1 & 1 \\ \hline
& 1 \\ \hline
\end{tabular}

\bigskip

\bigskip

(ii)\qquad $S_{4}\rightarrow S_{3}\rightarrow S_{2}$

\qquad \qquad
\begin{tabular}{|l|l|l|lllll|l|l|l|}
\cline{1-3}\cline{5-6}\cline{9-11}
1 &  &  &  & \multicolumn{1}{|l}{1} & \multicolumn{1}{|l}{} &
\multicolumn{1}{|l}{} &  &  & 1 &  \\ \cline{1-3}\cline{5-6}\cline{9-11}
1 & 1 &  &  & \multicolumn{1}{|l}{1} & \multicolumn{1}{|l}{1} &
\multicolumn{1}{|l}{} &  &  & 2 & 1 \\ \cline{1-3}\cline{5-6}\cline{9-11}
& 1 &  &  & \multicolumn{1}{|l}{} & \multicolumn{1}{|l}{1} &
\multicolumn{1}{|l}{} & = &  & 1 & 1 \\ \cline{1-3}\cline{5-6}\cline{9-11}
& 1 & 1 &  &  &  &  &  &  & 1 & 2 \\ \cline{1-3}\cline{9-11}
&  & 1 &  &  &  &  &  &  &  & 1 \\ \cline{1-3}\cline{9-11}
\end{tabular}

\bigskip

(iii)\qquad Finally

\qquad \qquad
\begin{tabular}{|l|l|l|l|l|l|l|l|lll|l|l|}
\cline{1-5}\cline{7-8}\cline{12-13}
1 &  &  &  &  &  & 1 &  &  &  &  & 1 &  \\
\cline{1-5}\cline{7-8}\cline{12-13}
&  & 1 & 1 &  &  &  & 1 &  &  &  & 2 & 1 \\
\cline{1-5}\cline{7-8}\cline{12-13}
1 & 1 &  &  &  &  & 1 &  &  & = &  & 1 & 1 \\
\cline{1-5}\cline{7-8}\cline{12-13}
&  &  & 1 & 1 &  & 1 & 1 &  &  &  & 1 & 2 \\
\cline{1-5}\cline{7-8}\cline{12-13}
& 1 &  &  &  &  &  & 1 &  &  &  &  & 1 \\ \cline{1-5}\cline{7-8}\cline{12-13}
\end{tabular}

\bigskip

Let us remark that Method (B) can be employed to verity the branching result
obtained by following Method (A).

\section*{Acknowledgments}

This work was supported in part by Fondecyt (Project 1160305).

\newpage

\begin{flushleft}
{\Large \textbf{Appendix (I)}}
\end{flushleft}

\begin{itemize}
\item[(\textit{i})] Let $g$ be an element of $S_{2N}$ and $C\left( g\right) $
the conjugacy class of $g$ in $S_{2N}$. The character $F_{W\left(
B_{N}\right) }^{S_{2N}}$ may be defined as follows:%
\begin{equation*}
F_{W\left( B_{N}\right) }^{S_{2N}}=\frac{\left\vert S_{2N}\right\vert }{%
\left\vert W\left( B_{N}\right) \right\vert }\frac{\left\vert C\left(
g\right) \cap W\left( B_{N}\right) \right\vert }{\left\vert C\left( g\right)
\right\vert }
\end{equation*}%
where $\left\vert S_{2N}\right\vert $ and $\left\vert W\left( B_{N}\right)
\right\vert $ are the orders of $S_{2N}$ and $W\left( B_{N}\right) $ and $%
\left\vert C\left( g\right) \cap W\left( B_{N}\right) \right\vert $ and $%
\left\vert C\left( g\right) \right\vert $ are, respectively, the orders of
the class $g$ in $W\left( B_{N}\right) $ and the order of the class $g$ of $%
S_{2N}\,$. Hence%
\begin{equation*}
F_{W\left( B_{N}\right) }^{S_{2N}}=\frac{\left( 2N\right) !}{2^{N}N!}\frac{%
\left\vert C\left( g\right) \cap W\left( B_{N}\right) \right\vert }{%
\left\vert C\left( g\right) \right\vert }\,.
\end{equation*}

\item[(\textit{ii})] For an even number $2N$ the number of partitions whose
subpartitions are even numbers is $P\left( N\right) $. For instance for $N=4$%
,%
\begin{equation*}
P\left( 8\right) =\left( 8\right) +\left( 6,2\right) +\left( 4,4\right)
+\left( 4,2,2\right) +\left( 2,2,2,2\right) =5=P\left( 4\right) \,.
\end{equation*}

\item[(\textit{iii})] Let the irreducible characters of $S_{2N}$
corresponding to such partitions compose $F_{W\left( B_{N}\right) }^{S_{2N}}$%
. For $N=4,$ \ \ $F_{4}=x\left( 4\right) +x\left( 2,2\right) $. By means of
the irreducible character table of $S_{4}$ it is possible to write:

\item[ ]
\begin{tabular}{cccccc}
order &  &  & $x(4)$ & $x(2,2)$ & $F_{W\left( B_{2}\right) }^{S_{4}}$ \\
\cline{2-6}
1 & \multicolumn{1}{|c}{$C_{1}$} & \multicolumn{1}{|c}{$1^{4}$} &
\multicolumn{1}{|c}{1} & \multicolumn{1}{|c}{2} & \multicolumn{1}{|c|}{3} \\
\cline{2-6}
6 & \multicolumn{1}{|c}{$C_{2}$} & \multicolumn{1}{|c}{$1^{2}2$} &
\multicolumn{1}{|c}{1} & \multicolumn{1}{|c}{0} & \multicolumn{1}{|c|}{1} \\
\cline{2-6}
3 & \multicolumn{1}{|c}{$C_{3}$} & \multicolumn{1}{|c}{$2^{2}$} &
\multicolumn{1}{|c}{1} & \multicolumn{1}{|c}{2} & \multicolumn{1}{|c|}{3} \\
\cline{2-6}
8 & \multicolumn{1}{|c}{$C_{4}$} & \multicolumn{1}{|c}{$13$} &
\multicolumn{1}{|c}{1} & \multicolumn{1}{|c}{-1} & \multicolumn{1}{|c|}{0}
\\ \cline{2-6}
6 & \multicolumn{1}{|c}{$C_{5}$} & \multicolumn{1}{|c}{$4$} &
\multicolumn{1}{|c}{1} & \multicolumn{1}{|c}{0} & \multicolumn{1}{|c|}{1} \\
\cline{2-6}
\end{tabular}
\end{itemize}

\quad From the formulas stated in (i), $\left\vert C\left( g\right) \cap
W\left( B_{2}\right) \right\vert $ can be evaluated:

\quad order of $C_{1}$ in $W\left( B_{2}\right) $ = 1

\quad order of $C_{2}$ in $W\left( B_{2}\right) $ = 2

\quad order of $C_{3}$ in $W\left( B_{2}\right) $ = 3

\quad order of $C_{4}$ in $W\left( B_{2}\right) $ = 0

\quad order of $C_{5}$ in $W\left( B_{2}\right) $ = 2

\begin{itemize}
\item[ ] The order of $W\left( B_{2}\right) $ is 2$^{2}$2!=8. The order of $%
C_{3}$ does not divide the order of $W\left( B_{2}\right) $. So the class $%
C_{3}$ of $W\left( B_{2}\right) $ must be decomposed in the character table
of $S_{4}$ and the class $C_{4}$ must be omitted. The resulting irreducible
character table of $S_{4}$ $\left( \text{denoted }X^{\prime }\right) $ is:

\item[ ]
\begin{tabular}{lll|c|c|c|c|c|}
\cline{4-8}
&  &  & 1$^{4}$ & 1$^{2}$2 & 2$^{2}$ & 2$^{2}$ & 4 \\ \cline{4-8}
&  &  & 1 & 1 & 1 & 1 & 1 \\ \cline{4-8}
&  &  & 3 & 1 & -1 & -1 & -1 \\ \cline{4-8}
$X\prime $ & = &  & 2 & 0 & 2 & 2 & 0 \\ \cline{4-8}
&  &  & 3 & -1 & -1 & -1 & 1 \\ \cline{4-8}
&  &  & 1 & -1 & 1 & 1 & -1 \\ \cline{4-8}
\end{tabular}
\end{itemize}

Remarks (1): For the identity class 1$^{2N}$ the character $F_{W\left(
B_{N}\right) }^{S_{2N}}$ is:

$N=2\qquad \qquad F_{W\left( B_{2}\right) }^{S_{4}}=3=3\cdot 1$

$N=3\qquad \qquad F_{W\left( B_{3}\right) }^{S_{6}}=15=5\cdot 3\cdot 1$

$N=4\qquad \qquad F_{W\left( B_{4}\right) }^{S_{8}}=105=7\cdot 5\cdot 3\cdot
1$

$N=5\qquad \qquad F_{W\left( B_{5}\right) }^{S_{10}}=945=9\cdot 7\cdot
5\cdot 3\cdot 1$

Accordingly:

\qquad \qquad $F_{W\left( B_{N}\right) }^{S_{2N}}=\left( 2N-1\right) \left(
2N-3\right) \ldots 1$

\bigskip

(2): In general if $\left\vert C\left( g\right) \cap W\left( B_{N}\right)
\right\vert $ is not a divisor of $\left\vert W\left( B_{N}\right)
\right\vert $ the corresponding class in $X^{\prime }\left( S_{2N}\right) $
must be divided.

For $N=3$ this occurs for the classes $\left( 1^{2}2^{2}\right) $ and $%
\left( 2^{3}\right) $; for $N=4$, the classes $\left( 1^{4}2^{2}\right) ,$ $%
\left( 1^{2}2^{3}\right) ,$ $\left( 2^{2}4\right) $ and $\left( 4^{2}\right)
$ are decomposed. It must be emphazided that for each $N$ the procedure must
be carried out. Perhaps this is the main difficulty of the present algorithm
for the reduction $S_{2N}\rightarrow W\left( B_{N}\right) .$

(3): The induced character table of $S_{2N}$, $\phi $, is treated in an
analogous manner. A modified character table, $\phi ^{\prime }$, results. So
for $S_{4}$:

\bigskip

\begin{tabular}{lllccccc}
&  &  & 1$^{4}$ & 1$^{2}$2 & 2$^{2}$ & 2$^{2}$ & 1$^{4}$ \\ \cline{4-8}
&  &  & \multicolumn{1}{|c}{1} & \multicolumn{1}{|c}{1} &
\multicolumn{1}{|c}{1} & \multicolumn{1}{|c}{1} & \multicolumn{1}{|c|}{1} \\
\cline{4-8}
&  &  & \multicolumn{1}{|c}{4} & \multicolumn{1}{|c}{2} &
\multicolumn{1}{|c}{0} & \multicolumn{1}{|c}{0} & \multicolumn{1}{|c|}{0} \\
\cline{4-8}
$\phi ^{\prime }$ & = &  & \multicolumn{1}{|c}{6} & \multicolumn{1}{|c}{2} &
\multicolumn{1}{|c}{2} & \multicolumn{1}{|c}{2} & \multicolumn{1}{|c|}{0} \\
\cline{4-8}
&  &  & \multicolumn{1}{|c}{12} & \multicolumn{1}{|c}{2} &
\multicolumn{1}{|c}{0} & \multicolumn{1}{|c}{0} & \multicolumn{1}{|c|}{0} \\
\cline{4-8}
&  &  & \multicolumn{1}{|c}{24} & \multicolumn{1}{|c}{0} &
\multicolumn{1}{|c}{0} & \multicolumn{1}{|c}{0} & \multicolumn{1}{|c|}{0} \\
\cline{4-8}
\end{tabular}

\bigskip

$\phi ^{\prime }$ and $X^{\prime }$ are related by the equation:

$\phi ^{\prime }=\Delta ^{\prime }X^{\prime }$\bigskip

\begin{tabular}{|l|l|l|l|l|lll|l|l|l|l|l|l|l|l|l|l|l|}
\cline{1-5}\cline{9-13}\cline{15-19}
1 & 1 & 1 & 1 & 1 &  &  &  & 1 & 0 & 0 & 0 & 0 &  & 1 & 1 & 1 & 1 & 1 \\
\cline{1-5}\cline{9-13}\cline{15-19}
4 & 2 & 0 & 0 & 0 &  &  &  & 1 & 1 & 0 & 0 & 0 &  & 3 & 1 & -1 & -1 & -1 \\
\cline{1-5}\cline{9-13}\cline{15-19}
6 & 2 & 2 & 2 & 0 &  & = &  & 1 & 1 & 1 & 0 & 0 &  & 2 & 0 & 2 & 2 & 0 \\
\cline{1-5}\cline{9-13}\cline{15-19}
12 & 2 & 0 & 0 & 0 &  &  &  & 1 & 2 & 1 & 1 &  &  & 3 & -1 & -1 & 1 & 1 \\
\cline{1-5}\cline{9-13}\cline{15-19}
24 & 0 & 0 & 0 & 0 &  &  &  & 1 & 3 & 2 & 3 & 1 &  & 1 & -1 & 1 & 1 & -1 \\
\cline{1-5}\cline{9-13}\cline{15-19}
\end{tabular}

\bigskip

Note that $\Delta =\Delta ^{\prime }$.

\newpage

\begin{flushleft}
\bigskip {\Large \textbf{Appendix (II)}}

The Reduction $S_{6}\rightarrow W\left( B_{3}\right) $ (Method (B))
\end{flushleft}

(1)$\qquad W\left( B_{3}\right) \rightarrow W\left( B_{2}\right) \rightarrow
W\left( B_{1}\right) $

\begin{tabular}{|l|l|l|l|l|llllll|l|l|}
\cline{1-5}\cline{7-8}\cline{12-13}
1 &  &  &  &  &  & \multicolumn{1}{|l}{1} & \multicolumn{1}{|l}{} &
\multicolumn{1}{|l}{} &  &  & 1 &  \\ \cline{1-5}\cline{7-8}\cline{12-13}
& 1 &  &  &  &  & \multicolumn{1}{|l}{} & \multicolumn{1}{|l}{1} &
\multicolumn{1}{|l}{} &  &  &  & 1 \\ \cline{1-5}\cline{7-8}\cline{12-13}
1 &  & 1 &  &  &  & \multicolumn{1}{|l}{1} & \multicolumn{1}{|l}{} &
\multicolumn{1}{|l}{} & = &  & 2 &  \\ \cline{1-5}\cline{7-8}\cline{12-13}
1 &  &  & 1 &  &  & \multicolumn{1}{|l}{1} & \multicolumn{1}{|l}{1} &
\multicolumn{1}{|l}{} &  &  & 2 & 1 \\ \cline{1-5}\cline{7-8}\cline{12-13}
& 1 &  & 1 &  &  & \multicolumn{1}{|l}{} & \multicolumn{1}{|l}{1} &
\multicolumn{1}{|l}{} &  &  & 1 & 2 \\ \cline{1-5}\cline{7-8}\cline{12-13}
& 1 &  &  & 1 &  &  &  &  &  &  &  & 2 \\ \cline{1-5}\cline{12-13}
&  & 1 &  &  &  &  &  &  &  &  & 1 &  \\ \cline{1-5}\cline{12-13}
&  & 1 & 1 &  &  &  &  &  &  &  & 2 & 1 \\ \cline{1-5}\cline{12-13}
&  &  & 1 & 1 &  &  &  &  &  &  & 1 & 2 \\ \cline{1-5}\cline{12-13}
&  &  &  & 1 &  &  &  &  &  &  &  & 1 \\ \cline{1-5}\cline{12-13}
\end{tabular}

\bigskip

(2)$\qquad S_{6}\rightarrow S_{2}$

\begin{tabular}{|l|l|l|llllll|l|l|}
\cline{1-3}\cline{5-6}\cline{10-11}
1 &  &  &  & \multicolumn{1}{|l}{1} & \multicolumn{1}{|l}{} &
\multicolumn{1}{|l}{} &  &  & 1 &  \\ \cline{1-3}\cline{5-6}\cline{10-11}
3 & 1 &  &  & \multicolumn{1}{|l}{1} & \multicolumn{1}{|l}{1} &
\multicolumn{1}{|l}{} &  &  & 4 & 1 \\ \cline{1-3}\cline{5-6}\cline{10-11}
3 & 3 &  &  & \multicolumn{1}{|l}{} & \multicolumn{1}{|l}{1} &
\multicolumn{1}{|l}{} &  &  & 6 & 3 \\ \cline{1-3}\cline{5-6}\cline{10-11}
3 & 3 & 1 &  &  &  &  &  &  & 6 & 4 \\ \cline{1-3}\cline{10-11}
1 & 2 &  &  &  &  &  &  &  & 3 & 2 \\ \cline{1-3}\cline{10-11}
2 & 6 & 2 &  &  &  &  & = &  & 8 & 8 \\ \cline{1-3}\cline{10-11}
1 & 3 & 3 &  &  &  &  &  &  & 4 & 6 \\ \cline{1-3}\cline{10-11}
& 2 & 1 &  &  &  &  &  &  & 2 & 3 \\ \cline{1-3}\cline{10-11}
& 3 & 3 &  &  &  &  &  &  & 3 & 6 \\ \cline{1-3}\cline{10-11}
& 1 & 2 &  &  &  &  &  &  & 1 & 4 \\ \cline{1-3}\cline{10-11}
&  & 1 &  &  &  &  &  &  &  & 1 \\ \cline{1-3}\cline{10-11}
\end{tabular}

\bigskip \newpage

(3)$\qquad \left\{ S_{6}\rightarrow W\left( B_{3}\right) \right\} \left\{
W\left( B_{3}\right) \rightarrow W\left( B_{1}\right) \right\} =\left\{
S_{6}\rightarrow S_{2}\right\} $

\begin{tabular}{|l|l|l|l|l|l|l|l|l|l|llllll|l|l|}
\cline{1-10}\cline{12-13}\cline{17-18}
1 &  &  &  &  &  &  &  &  &  &  & \multicolumn{1}{|l}{1} &
\multicolumn{1}{|l}{} & \multicolumn{1}{|l}{} &  &  & 1 &  \\
\cline{1-10}\cline{12-13}\cline{17-18}
&  & 1 & 1 &  &  &  &  &  &  &  & \multicolumn{1}{|l}{} &
\multicolumn{1}{|l}{1} & \multicolumn{1}{|l}{} &  &  & 4 & 1 \\
\cline{1-10}\cline{12-13}\cline{17-18}
1 &  & 1 &  & 1 &  &  & 1 &  &  &  & \multicolumn{1}{|l}{2} &
\multicolumn{1}{|l}{} & \multicolumn{1}{|l}{} &  &  & 6 & 3 \\
\cline{1-10}\cline{12-13}\cline{17-18}
&  &  & 1 &  &  & 1 & 1 & 1 &  &  & \multicolumn{1}{|l}{2} &
\multicolumn{1}{|l}{1} & \multicolumn{1}{|l}{} &  &  & 6 & 4 \\
\cline{1-10}\cline{12-13}\cline{17-18}
& 1 &  & 1 &  &  & 1 &  &  &  &  & \multicolumn{1}{|l}{1} &
\multicolumn{1}{|l}{2} & \multicolumn{1}{|l}{} &  &  & 3 & 2 \\
\cline{1-10}\cline{12-13}\cline{17-18}
&  & 1 & 1 & 1 & 1 &  & 1 & 1 &  &  & \multicolumn{1}{|l}{} &
\multicolumn{1}{|l}{2} & \multicolumn{1}{|l}{} & = &  & 8 & 8 \\
\cline{1-10}\cline{12-13}\cline{17-18}
&  &  &  & 1 &  &  & 1 & 1 & 1 &  & \multicolumn{1}{|l}{1} &
\multicolumn{1}{|l}{} & \multicolumn{1}{|l}{} &  &  & 4 & 6 \\
\cline{1-10}\cline{12-13}\cline{17-18}
1 &  &  &  & 1 &  &  &  &  & 1 &  & \multicolumn{1}{|l}{2} &
\multicolumn{1}{|l}{1} & \multicolumn{1}{|l}{} &  &  & 2 & 3 \\
\cline{1-10}\cline{12-13}\cline{17-18}
& 1 &  & 1 &  & 1 &  &  & 1 &  &  & \multicolumn{1}{|l}{1} &
\multicolumn{1}{|l}{2} & \multicolumn{1}{|l}{} &  &  & 3 & 6 \\
\cline{1-10}\cline{12-13}\cline{17-18}
&  &  &  & 1 & 1 &  &  &  &  &  & \multicolumn{1}{|l}{} &
\multicolumn{1}{|l}{1} & \multicolumn{1}{|l}{} &  &  & 1 & 4 \\
\cline{1-10}\cline{12-13}\cline{17-18}
& 1 &  &  &  &  &  &  &  &  &  &  &  &  &  &  &  & 1 \\
\cline{1-10}\cline{17-18}
\end{tabular}

\end{document}